\documentclass[10pt,technote]{IEEEtran}
\usepackage{amsmath,amssymb,amsfonts}
\usepackage{color,graphicx}
\usepackage[version=3]{mhchem}
\usepackage{paralist}
\newcommand{\reffig}[1]{Figure \ref{#1}}
\newcommand{\ImagenA}[5]{\begin{figure*}[#4]
  \centering
  \includegraphics[angle=0,width=#5\textwidth]{#2}
  \caption{#1}
  \label{#3}
  \end{figure*}}

\newcommand{\Reaction}[3]{#1 \ce{->[#3]} #2}
\newcommand{\ReactionR}[4]{#1 \ce{<=>[#3][#4]} #2}

\newcommand{\bs}[1]{\boldsymbol {#1}}

\newcommand{\con}{\ensuremath{c}}

\newcommand{\Sb}{\ensuremath{S}}
\newcommand{\sub}{\ensuremath{s}}
\newcommand{\Sub}{\ensuremath{\mathbf{\sub}}}
\newcommand{\SUB}{\ensuremath{\mathbf{S}}}
\newcommand{\proba}{\ensuremath{p}}
\newcommand{\Prob}{\ensuremath{\mathbf{\proba}}}

\newcommand{\kf}[1]{\ensuremath{k_{\mathrm{f} #1}}}
\newcommand{\kb}[1]{\ensuremath{k_{\mathrm{b} #1}}}

\newcommand{\te}{\ensuremath{\mathnormal{t}}}

\newtheorem{Theorem}{Theorem}[section]

\newcommand{\dt}{\mathrm{d}\te}

\newcommand{\ddt}{\frac{\mathrm{d}}{\dt}}

\newcommand{\iga}{&=}
\newcommand{\0}{\ensuremath{\bs 0}}

\newcommand{\abs}[1]{\ensuremath{\left| #1 \right| }}

\usepackage[colorinlistoftodos,textsize=small,textwidth=2cm]{todonotes}

\newcommand{\trace}{\ensuremath{\mathrm{trace}}}

\newcommand{\OpPr}[1]{\ensuremath{\mathrm{Pr\left ( #1\right )}}}

\newcommand{\orders}{n}
\newcommand{\nums}{w}
\newcommand{\thh}{\ensuremath{^{th}}~}

\graphicspath{ {./Figures/} }

\begin{document}
\title{Order Reduction of the Chemical Master Equation via Balanced Realisation} 
\author{Fernando L\'opez-Caamal, Tatiana T. Marquez-Lago$^*$\\ Integrative Systems Biology Unit, Okinawa Institute of Science and Technology, Kunigami, Okinawa 904-0412, Japan\\
$^*$ Author for correspondence: tatiana.marquez@oist.jp}
\date{\today}

\maketitle 

\begin{abstract}
We consider a Markov process in continuous time with a finite number of discrete states. The time-dependent probabilities of being in any state of the Markov chain are governed by a set of ordinary differential equations, whose dimension might be large even for trivial systems. Here, we derive a reduced ODE set that accurately approximates the probabilities of subspaces of interest with a known error bound. Our methodology is based on model reduction by balanced truncation and can be considerably more computationally efficient than the Finite State Projection Algorithm (FSP)
 when used for obtaining transient responses. We show the applicability of our method by analysing stochastic chemical reactions. First, we obtain a reduced order model for the infinitesimal generator of a Markov chain that models a reversible, monomolecular reaction. In such an example, we obtain an approximation of the output of a model with $301$ states by a reduced model with $10$ states. Later, we obtain a reduced order model for a catalytic conversion of substrate to a product; and compare its dynamics with a stochastic Michaelis-Menten representation. For this example, we highlight the savings on the computational load obtained by means of the reduced-order model. Finally, we revisit the substrate catalytic conversion by obtaining a lower-order model that approximates the probability of having predefined ranges of product molecules.
\end{abstract}

\section{Introduction}
Markov chains are versatile dynamical systems that model a broad spectrum of physical, biological, and engineering systems. Along with the broad range of its applications, one of the main advantages of Markov chains is that some of them can be easily handled and cast as time-invariant, linear systems \cite{feller2008introduction,pierre1999markov,van2007stochastic}. 

In this paper, we focus on continuous-time, discrete-state, homogeneous, irreducible Markov chains with a finite number of states. The probability of being in any state is governed by a linear set of ordinary differential equations (ODEs), where individual ODEs correspond to each state of the system, describing all possible transitions in and out of such states. This set of ODEs is commonly referred to as forward Kolmogorov equation or Chemical Master equation and might have large dimensions even for simple systems. Hence, obtaining a solution for such a system might be analytically intractable and computationally demanding.  

Provided that one is interested only in some states or a combination of states of the Markov chain, it is possible to obtain a reduced order model via the balanced realisation of the linear system that describes the probability of being in such states. The reduced model has a smaller number of coupled differential equations, yet approximates the output of the full model with an error bound proportional to the sum of the Hankel Singular Values neglected to obtain the reduced model \cite{Moore1981,Zhou1996robust,Gugercin2004,Skogestad2007}. Since chemical reaction networks in a homogeneous media with a low number of molecules, and in thermodynamic equilibrium can be described as Markov chains, it is possible to apply our methodology to this class of systems.

There exist alternative approaches to obtain reduced order models from the CME. For instance, the Finite State Projection Method obtains the probability density function in prescribed subsets of the state space and for a specific time point \cite{Munsky2006}. From this method, it is also possible to obtain a linear set of ODEs that can, in general, be further reduced via the methodology we use in this paper. Other approaches make use of Krylov subspaces to approximate the solution of the exponential matrix that generates the solutions of the Markov chain \cite{saad1992analysis, burrage2006krylov}. Additionally, when the species can be classified by its behaviour into stochastic or deterministic, \cite{menz2012hybrid} propose a methodology in which the CME can be solved directly and efficiently, when the number of species with stochastic behaviour is low. In this direction, works like \cite{haseltine2002approximate} avail of a time scale separation to estimate the solution of the fast-varying species; and use this estimation to approximate the trajectories of the slow-varying species. On a different perspective, \cite{jahnke2011reduced} analysed methods to approximate the solution of selected states of the CME, when such solutions can be expressed as the product of two probability density functions: one that describes probabilities of states of interest and a second that depends on the rest of states. This later probability distribution can be approximated by its mean, for instance, so as to yield an approximated probability density function for those probabilities of interest. However, this approach might yield coarse results if the underlying assumptions are crude. 

As an alternative, when the analytical or computational treatment of the Markov chain is infeasible, it is common to opt for numerical simulations of the stochastic system and analyse the outcome statistically. \cite{Barrio2010}, and \cite{ullah2011stochastic}, among many others, provide surveys of simulation methods of stochastic reaction networks. However, these methods might require large computational times to yield accurate results.

A different way to reduce the CMEs is to consider subsystems that focus on features of interest. In the chemical context, \cite{gillespie2009subtle} showed as a proof-of-concept that the simple reaction $\Reaction{\ReactionR{\Sb_{1}}{\Sb_2}{k_1}{k_2}}{\Sb_3}{k_3}$ can only be accurately represented by $\Reaction{\Sb_{1}}{\Sb_3}{k}$ under special conditions on the parameters $k_1,k_2,k_3$, which render the dynamics of the species $\Sb_2$ irrelevant for the behaviour of $\Sb_3$. This study highlights the shortcomings of neglecting species within a stochastic reaction network. In this paper, we adopt a different approach and overcome these difficulties by deriving a reduced-order model. Such reduced-order model accurately approximates the dynamics of the underlying Markov chain for selected states of any chain with any kind of reaction propensities. 

There exist, however, exact approaches that abridge specific topologies of reaction networks. For instance, in \cite{Barrios2013,leier2013sub} different classes of monomolecular reaction networks are exactly represented as reactions characterised by delay distributions. In turn, works like \cite{lopez2013exact,lee2012analytical} are committed to obtain exact analytical solutions of stochastic chemical reaction networks with linear and nonlinear reactions. Importantly, once a reduced ODE set via balanced realisation is obtained, one can avail of the results in \cite{lopez2013exact} to derive a closed-form expression for approximation of the CME solution. 

We illustrate our methodology with the analysis of a reversible, stochastic reaction whose CME has $301$ states. In contrast, an adequate reduced order model has only $10$ states and yield an $\mathcal L_2$ gain of the approximation error of $587.9172\times 10^{-6}$. Later, we obtain a reduced order model that approximates the catalysed conversion of a substrate to a product, even in cases in which a stochastic Michaelis-Menten approximation fails to obtain accurate results \cite{sanft2011legitimacy}. For such a system, 
the simulation of the reduced model may be several orders of magnitude faster than the simulation of the CME. However, there exist an initial cost in computational time to derive the reduced order model. Hence, obtaining the model reduction is profitable when the lower-dimensional ODE set is used repeatedly. Finally, we derive a model that approximates the probability of having predefined ranges of product molecules, in the same catalytic substrate conversion.
\section{Mathematial Background}\label{Sec:Back}
\subsection{Continuous-Time, Discrete-State, Homogeneous Markov Chain}\label{Sec:MarCh}
Consider a discrete and finite set of states 
\begin{align}\label{Eq:SetSpace}
\SUB = \left \{ \Sub^{i} \in \mathbb Z^\orders, \: \forall \, i \in [1,\nums] \right \}
\end{align}
and let the system's state at time $\te$ be denoted by $\Sub(\te) : \mathbb R_+ \to \SUB$. Moreover, we consider that the transition from one state to another can be modelled by a time-homogeneous Markov chain, i.e., the next state, $\Sub(\te+\dt)$, only depends on the current state, \Sub(\te), independently of $\te$. We use $\proba_i(\te):\mathbb R_+ \times \SUB \to [0,1]\subset \mathbb R$ to denote the probability, $\OpPr{\circ}$, of the system's state to be $\Sub^i$ at time $\te$. This notation and the Markov property add up to 
\begin{equation*}
 	\proba_i(\te+\dt) = \mathrm{Pr} \left ( \Sub(\te+\dt) = \Sub^i | \Sub(\te)  \right).
\end{equation*}
We gather the probabilities for every state in the column vector 
\begin{equation}\label{Eq:VecP}
\Prob(\te) := (\proba_1(\te) \quad \hdots \quad \proba_\nums(\te) )^T.
\end{equation}
Let us denote the transition probability from state $j$ to state $i$ at time $\te+\tau$ by $q_{ij}(\te+\tau):\mathbb R_+ \times \SUB^2\to [0,1]$. That is to say,
\begin{equation*}
	q_{ij}(\te+\tau):=\OpPr{\Sub(\te+\tau) = \Sub^i | \Sub(\te) = \Sub^j}.
\end{equation*}
The time-homogeneity property of the Markov chain implies
 \begin{align}\label{Eq:tranij}
 	q_{ij}(\te+\tau)&:=\nonumber\\
 	\OpPr{\Sub(\te+\tau) = \Sub^i | \Sub(\te) = \Sub^j}&=\nonumber\\
    \sum_{k=1}^{\nums} \OpPr{\Sub(\te+\tau) = \Sub^i; \Sub(\tau) = \Sub^k | \Sub(0) = \Sub^j}&=\nonumber \\
 							\sum_{k=1}^{\nums} \OpPr{\Sub(\te+\tau) = \Sub^i | \Sub(\tau) = \Sub^k ; \Sub(0) = \Sub^j}\nonumber \\  \times \OpPr{\Sub(\tau)=\Sub^k | \Sub(0) = \Sub^j}&=\nonumber\\
 							\sum_{k=1}^{\nums} q_{ik}(\te)q_{kj}(\tau)&.
\end{align}
In matrix form \eqref{Eq:tranij}, known as \emph{Chapman-Kolmogorov equation}, is
\begin{align}
\mathbb R_+ \to [0,1]^{\nums\times\nums} : \mathbf{Q}(\te+\tau) = \mathbf{Q}(\te) \mathbf {Q}(\tau)  =\mathbf{Q}(\tau) \mathbf {Q}(\te). \label{Eq:SemiGroup}
\end{align}
This matrix gathers all the transition probabilities as a function of time and, by consequence, its columns add to one for all $\te$. Additionally, if the Markov chain is \emph{irreducible}, $\mathbf Q$ has a simple eigenvalue  $\lambda_1 = 1$, and $\lambda_1 > \abs{\lambda_i}\, \forall\,  1 <  i  \leq w$. This is consequence of the Perron-Frobenius Theorem as described in \cite[Ch. 6]{pierre1999markov}, for example. In the rest of this paper, we will deal with finite, irreducible, homogeneous, continuous-time, discrete-state Markov chains exclusively. 

Our main interest is to determine the time-dependent probabilities of being in any state of the chain. To this end, we consider the  \emph{infinitesimal generator} of the Markov chain defined as
\begin{align} \label{Eq:InfGen}
	\mathcal{A} :=& \lim_{\tau \to 0} \frac{\mathbf Q(\tau) - \mathbf{I}}{\tau}.
\end{align}
The elements of the matrix above are
\begin{align}\label{Eq:InfGenEW}
	a_{ij} \iga \lim_{\tau \to 0} \frac{ q_{ij}(\tau) - \delta_{ij} }{\tau},
\end{align}
where it can be shown that the elements $a_{ij}$ satisfy
\begin{align*}
	a_{ii} = - \sum_{j = 1, j \neq i}^{\nums} a_{ij}.
\end{align*} 
The last relationship above shows that every column of  $\mathcal A$ adds up to zero, provided each column of $\mathbf{Q}(\circ)$ add up to one. 


It is well-known that $\mathcal A$ is the generator of the positive semigroup that governs the evolution of \Prob(\te) (see \cite[Sec. 5.6]{allen2003introduction}, for instance):
\begin{align}\label{Eq:ODEp}
 	\ddt \Prob(\te) = \mathcal{A} \Prob(\te) , \quad \Prob(0) = \Prob_0.
\end{align}
Under our assumptions, the Markov chain is irreducible and with a finite number of states. Hence, $\mathbf Q$ has a unique Frobenius eigenvalue with algebraic multiplicity one \cite{farina2011positive}. The simple Perron-Frobenius eigenvalue of the stochastic matrix $\mathbf Q$  in \eqref{Eq:SemiGroup} is 1. Now, let $\bs \nu$ and $\lambda$ be the right Perron-Frobenius eigenvector and eigenvalues of $\mathbf Q$, then the eigenvalues of $\mathcal A$ satisfy
\begin{align*}
	\mathcal A \bs \nu \iga ( \mathbf Q - \mathbf I ) \bs \nu / \tau,\\
	 								 \iga ( \lambda - 1 ) \bs \nu / \tau.
\end{align*}
That is, $\mathcal A$ preserves the configuration of the eigenvalues of $\mathbf{Q}$, upon shifting one unit to the left and rescaling. This implies that $\mathcal A$ has a zero eigenvalue and the rest of its eigenvalues have negative real part, as confirmed by analysing the Gershgorin circles of the columns of $\mathcal A$. 

Note that the dimension of $\Prob(\circ)$, $\nums$, might be large as it represents all the configurations of a system with $\orders$ characteristics. In the population and biochemical context, $\orders$ represents the number of species, whereas $\nums$ is the number of all the possible combination of species' population counts. In the following section, we model a stochastic chemical reaction network with the Markov chains described above.

\subsection{Chemical Master Equation}
Now, let us consider $\orders$ species in a homogeneous medium and in thermodynamic equilibrium and a set of $m$ reactions represented by 
\begin{equation}
 \Reaction{\sum_{i=1}^\orders \alpha_{ij} S_i}{\sum_{i=1}^\orders \beta_{ij} S_i}{a_j}.\label{Eq:Reac}
 \end{equation}
Let the entries of the stoichiometric matrix $\mathbf N \in \mathbb N^{\orders \times m}$ be 
\begin{align*}
 	n_{ij} := \beta_{ij} - \alpha_{ij}. 
\end{align*}
Furthermore, let us consider a vector comprised of the number of molecules, $\sub_i(\te)$, for every species, $\Sb_i$:
\begin{equation}
\Sub(\te) = \left ( \sub_1(\te)\quad\sub_2(\te)\quad\hdots\quad\sub_\orders (\te)  \right )^{T} : \mathbb R_+ \to \SUB.
\end{equation}
The finite set $\SUB$ above was defined in \eqref{Eq:SetSpace} and contains, at least, all the possible combinations of the species' molecular numbers in the reaction network. Consider that the $i$\thh reaction is the only reaction happening within the interval $(\te , \te + \tau]$. Hence the number of molecules at time $t+\tau$ is
\begin{equation}
 \Sub (\te + \tau) = \Sub (\te) + \mathbf{n}^{i},
\end{equation}   
where $\mathbf{n}^{i}$ represents the $i$\thh column of $\mathbf{N}$. 

This reaction network may be modelled by the continuous-time, discrete-state jump Markov process described in Section \ref{Sec:MarCh}. The states of the Markov chain are the elements in \SUB. In turn, the vector $\Prob(\te)$ in \eqref{Eq:VecP} gathers the time-dependent probabilities of being in every state, whose time evolution is governed by \eqref{Eq:ODEp}. Additionally, the stochastic behaviour of thermally stable and spatially homogeneous reaction networks has been described\cite{Gillespie1992}, where the transition rates between states of the system are shown in Table \ref{Tb:Prope}. 
\begin{table}[b]
\caption{Reactions and their propensity function. The symbol $\sub_i$ denotes the number of molecules of the species $\Sb_i$. The symbol $a_i$ denotes the propensity of reaction $i$.} \label{Tb:Prope}
\begin{center}
\begin{tabular}{ll}
Reaction & Propensity \\
\hline
$\Reaction{0}{\Sb_1}{k_1}$ & $a_1 = k_1$\\ 
$\Reaction{\Sb_2}{X}{k_2}$ & $a_2(\sub_2(\te)) = k_2 \sub_2(\te) $ \\
$\Reaction{\Sb_3 + \Sb_4}{X}{k_3}$ & $a_3(\sub_3(\te),\sub_4(\te)) = k_3 \sub_3(\te) \sub_4(\te)$ \\
$\Reaction{2\Sb_5}{X}{k_4}$ & $a_4(\sub_5(\te)) =  k_4 \sub_5(\te) (\sub_5(\te)-1) \slash 2 \quad \sub_5(\te) \geq 1$ \\
\hline
\end{tabular}
\end{center}
\end{table}
To construct the matrix $\mathcal A$ in \eqref{Eq:ODEp}, we have to evaluate the probabilities' transition rate for all states $\mathbf \sub^i \in \SUB$ and arrange them as the entries of $\mathcal A$ as follows
\begin{equation}\label{Eq:PropMat}
\mathcal A_{ji} = \begin{cases} -\sum_{k=1}^{m} a_k(\Sub^i), & i=j, \\ a_k(\Sub^i),  & \forall\, j: \Sub^j = \Sub^i + \mathbf n^k,\\ 0, & \textnormal{otherwise}. \end{cases}
\end{equation} 
In the next section, we present a methodology used to obtain reduced order models, which are capable of reproducing the dynamical behaviour of a linear system with a smaller number of ODEs.

\subsection{Balanced Model Reduction}\label{Sec:Balanced}
In this section, we present an overview of a methodology used for obtaining lower dimensional models via balanced realisation. The literature on this topic is vast and we refer the interested reader to \cite{Skogestad2007,Zhou1996robust,Gugercin2004,benner2005model} for a deeper presentation of this type of model reduction.  

Let us consider a linear system of the form 
\begin{subequations}\label{Eq:LinODE}
\begin{align}
 \ddt \mathbf x (\te) \iga \mathbf{Ax (\te) +Bu (\te) },~\mathbf{x}(0) = \mathbf{x}_{0}\\  
 \mathbf{y}(\te) \iga \mathbf{Cx (\te) +Du (\te) }.
\end{align}
\end{subequations}
We will assume that all the eigenvalues of $\mathbf A$ have negative real part, i.e. stable, and that is both controllable and observable. \emph{Controllability} is the property of \eqref{Eq:LinODE} which ensures that it is possible to steer the state of the system from any initial condition $\mathbf{x}_{0}$ to any desired state at a specific time, by the application of an adequate forcing function $\mathbf u(\te)$. In turn, \emph{observability} refers to the capability of computing $\mathbf{x}_{0}$ given the knowledge of $\mathbf u(\te)$ and $\mathbf y(\te)$ for all previous time. These two properties hold when the (observability and controllability) matrices below are full rank\cite{chen1998linear}.
\begin{align*}
 	\mathcal O :=\begin{pmatrix}  \mathbf{C} \\  \mathbf{C} \mathbf{A} \\ \vdots \\  \mathbf{C} \mathbf{A}^{\nums-1} \end{pmatrix}, & \quad\quad \mathcal K := \begin{pmatrix}  \mathbf{B} & \mathbf{A}\mathbf{B} & \hdots & \mathbf{A}^{\nums-1} \mathbf{B}  \end{pmatrix} .
\end{align*}  
When \eqref{Eq:LinODE} is simultaneously stable, observable, and controllable there exist unique, symmetric, positive-definite matrices $\mathcal P$ and $\mathcal Q$ which are solution of the following Lyapunov equations
\begin{subequations}\label{Eq:Gram}
\begin{align}
 	\mathbf{A}\mathcal P  + \mathcal P \mathbf{A}^T + \mathbf{BB}^T \iga \bs 0,\\
 	\mathbf{A}^T \mathcal{Q} + \mathcal{Q} \mathbf{A}  + \mathbf{C}^T\mathbf C  \iga \bs 0.
\end{align}
\end{subequations}
%
The singular values of the product of $\mathcal P$ and $\mathcal Q$ are known as the Hankel singular values, $\sigma_i$, of the system. When 
\begin{align*}
	 \mathcal P =  \mathcal Q  = \mathrm{diag} (\sigma_1 ,\sigma_2 ,\hdots,\sigma_\nums),
\end{align*}
where $\sigma_1\leq\sigma_2\leq..\leq \sigma_\nums$, the system \eqref{Eq:LinODE} is Lyapunov balanced \cite{Moore1981}. 
If the original system \eqref{Eq:LinODE} is not stable nor represented by its minimal realisation (i.e.~ simultaneously observable and controllable), we suggest the transformation in \cite{Therapos1989} to obtain its Lyapunov balanced form.

An advantage of having a balanced realisation is that the magnitude of the singular values $\sigma_i$ decays quickly as $i$ increases. There are several techniques that avail of this observation to derive reduced order models, depending on the required characteristics of such a reduced model \cite{Gugercin2004}. Let us denote the coordinates of the balanced realisation of \eqref{Eq:LinODE} by $\mathbf{\tilde x}$ and, accordingly, we mark with a tilde the matrices related to the balanced realisation. One of the simplest approaches to obtain a reduced order model is to partition  $\mathbf{\tilde x}(\te)$ in \eqref{Eq:LinODE}, to obtain
\begin{align*}
\ddt \begin{pmatrix} \mathbf{\tilde x_1(\te)}\\\mathbf{\tilde x_2(\te)} \end{pmatrix} \iga  \begin{pmatrix} \mathbf{\tilde  A_{11}} & \mathbf{\tilde A_{12}} \\ \mathbf{\tilde A_{21}} & \mathbf{\tilde A_{22}}\end{pmatrix} \begin{pmatrix} \mathbf{\tilde x_1}(\te)\\\mathbf{\tilde x_2}(\te) \end{pmatrix} +\begin{pmatrix} \mathbf{\tilde B_1}\\\mathbf{\tilde B_2} \end{pmatrix} \mathbf{u(\te)},\\
\mathbf{y}(\te) \iga \begin{pmatrix} \mathbf{\tilde C_1}&\mathbf{\tilde C_2} \end{pmatrix}  \begin{pmatrix} \mathbf{\tilde x_1}(\te)\\\mathbf{\tilde x_2}(\te) \end{pmatrix} + \mathbf{\tilde D u(\te)}.
\end{align*}
This separation also induces the partition $\mathcal P = \mathcal Q = \mathrm{diag}(\Sigma_1,\Sigma_2)$, where $\Sigma_1 = \mathrm{diag}(\sigma_1,\hdots,\sigma_k)$ and $\Sigma_2 = \mathrm{diag}(\sigma_{k+1},\hdots,\sigma_{\nums})$. By neglecting the states associated to the small Hankel singular values, $\mathbf{\tilde x_2}$, the truncated model becomes \cite{Pernebo1982}
\begin{subequations}\label{Eq:ODEBalRed}
\begin{align}
 \ddt \mathbf{\tilde x_1}(\te) &\approx \mathbf{\tilde A_{11}}  \mathbf{\tilde x_1}(\te)  + \mathbf{\tilde B_1 u(\te)},\quad \mathbf{\tilde{x}_1}(0) =  \mathbf{\tilde{x}_{1\mathnormal{0}}} \\
 \mathbf{y}(\te) &\approx \mathbf{\tilde C_1} \mathbf{\tilde x_1}(\te) + \mathbf{\tilde Du(\te)}.
\end{align} 
\end{subequations}
This model is known to preserve the most important eigenvalues of the original system. However, some other properties such as steady state are slightly modified. When such a property is of interest, \emph{model reduction by residualisation} is more suitable \cite{Skogestad2007}. 
Both of these methods are already included in languages such as Phyton and Matlab, where the balanced realisation of a linear system is in the function \emph{balreal} and the model reduction via truncation and residualisation is in \emph{modred}. It is also important to mention, that the $\mathcal L_2$ gain of the approximation error for both model reduction methods above are as follows 
\begin{align}\label{Eq:ErrorBound}
 	\frac{\abs{\abs{ \mathbf y  -  \mathbf y_{\textrm{red}}}}_{\mathcal L_2}}{\abs{\abs{\mathbf u}}_{\mathcal L_2}} \leq 2  \sum_{i=k+1}^{\nums} \sigma_{i},
\end{align}
where $\mathbf y_{\textrm{red}}(\te)$ is the output of the reduced model. Please, refer to Appendix \ref{Sec:App} for a derivation of such an error bound. In the forthcoming section, we build upon the material in this section to obtain a reduced order model of the representation of a continuous-time, discrete-state, homogeneous, irreducible Markov chain. 
\section{Order Reduction of Infinitesimal Generators}\label{Sec:OrderRed}
In this section we are interested in the probability of being in some (linear combination of) states $\mathbf{y}(\te) : \mathbb R_+ \to \mathbb R^r \subset \mathbf S$ of the Markov chain. As noted in Equation \eqref{Eq:ODEp}, the vector $\Prob(\te)$ evolves according to the linear ODE
\begin{subequations}\label{Eq:ODEy}
\begin{align}
 	\ddt \Prob(\te) \iga \mathcal{A} \Prob(\te), \quad \Prob(0) = \Prob_0,\\
 	 \mathbf{y}(\te)  \iga  \mathcal{C} \Prob(\te). \label{Eq:Defye}
\end{align}
\end{subequations}
Also, as was mentioned in Section \ref{Sec:MarCh}, the infinitesimal generator $\mathcal A$ of an irreducible Markov chain with finite states has the properties :
\begin{subequations}\label{Eq:Properties}
\begin{align}
 \mathcal A &\leq 0 \quad(\textnormal{Its eigenvalues are nonpositive}), \\
\bs{1}^T \mathcal A   \iga \0^T \quad(\textnormal{Its colums add up to zero}).
\end{align}
\end{subequations}
Without loss of generality, we will assume that the system \eqref{Eq:ODEy} is both controllable and observable. Should it lack these two properties, there always exist a transformation that obtains the observable and controllable subspace of \eqref{Eq:ODEy}; namely, the Kalman Decomposition \cite{kalman1963mathematical, chen1998linear}.

To consider a reduced model that does not have a zero eigenvalue, let $\mathcal A$ be partitioned as follows:
\begin{align}\label{Eq:Apart}
\mathcal{A} = \left (\begin{array}{c|c}
													 a_{11} & \mathbf{a_{12}}^T\\\hline \mathbf{a_{21}} & \mathcal{A}_{22}										\end{array} 	\right ).
\end{align}
Also we consider the following similarity transformation
\begin{subequations}\label{Eq:Trans}
\begin{align}
 	 \mathbf{T} &:=
 														\left (\begin{array}{c|c}
													 1 & \bs{1}^T_{\nums -1}\\\hline \0_{\nums-1} & \mathbf I_{\nums-1 \times \nums -1} \end{array} 	\right )	\in \mathbb{N}^{\nums\times \nums},\\
 	\mathbf{T}^{-1} &:=
 														 \left (\begin{array}{c|c}
													 1 &  -\bs 1 ^T_{\nums -1} \\\hline \0_{\nums-1} & \mathbf I_{\nums-1 \times \nums -1} 	\end{array} 	\right ) \in \mathbb{N}^{\nums\times \nums}	,\\ 		
  \begin{pmatrix}z_0(\te)\\\mathbf{z}(\te) \end{pmatrix} &= \mathbf{T\Prob}(\te). 														 											 									\end{align} 
\end{subequations}
By differentiating the last equation above and using the expressions in \eqref{Eq:ODEy}, we get
\begin{align}
\ddt \begin{pmatrix} z_0(\te)\\\mathbf{z(\te)} \end{pmatrix} \iga  
\begin{pmatrix} 0 & \0^T \\ \mathbf{a_{21}} & \mathbf{A}\end{pmatrix} \begin{pmatrix} z_0(\te)\\\mathbf{z}(\te) \end{pmatrix} \\
z_0(0) \iga 1,~\mathbf z(0) = \begin{pmatrix}  \proba_2(0) & \hdots & \proba_\nums(0)  \end{pmatrix}^T,
 \label{Eq:ODEPaso}\\
\mathbf{y}(\te) \iga \mathcal{C}   \left (\begin{array}{c|c}
													 1 &  -\bs 1 ^T_{\nums -1} \\\hline \0_{\nums-1} & \mathbf I_{\nums-1 \times \nums -1}
													\end{array} 	\right )   \begin{pmatrix} z_0(\te)\\\mathbf{z}(\te) \end{pmatrix}\nonumber.
\end{align}
The solution for the first state is the unitary step function, that is $z_0(\te) = h(\te)$. By substituting this solution in the ODE above, we have
\begin{subequations}\label{Eq:ODEStable}
\begin{align}
 	\ddt \mathbf{z(\te)} \iga  \mathbf{A z}(\te) + \mathbf{b} h(\te),~\mathbf{z}(0) = \mathbf{z}_0 \\
\mathbf{y}(\te) \iga \mathbf{C}   \mathbf{z}(\te) + \mathbf{d} h(\te),
\end{align}
\end{subequations}
where
\begin{subequations}\label{Eq:ModStable}
\begin{align}
	\mathbf{A} &:= \mathcal{A}_{22} - \mathbf{a_{21}} \bs 1^T_{\nums-1},\\
	\mathbf{b} &:= \mathbf{a_{21}},\\
   \mathbf{C} &:= \mathcal{C} \begin{pmatrix} -\bs 1^T_{\nums-1}\\ \phantom{-}\mathbf{I}_{\nums-1}  \end{pmatrix},\\
   \mathbf{d} &:= \mathcal{C}\begin{pmatrix}  1\\\bs 0_{\nums -1} \end{pmatrix}.	 
\end{align}
\end{subequations}
The spectrum of $\mathbf A$ in \eqref{Eq:ODEStable} has all the eigenvalues of $\mathcal A$, except for the zero eigenvalue. To see this, recall that the trace of a matrix is the sum of its eigenvalues. As \eqref{Eq:ODEPaso} arises from a similarity transformation applied to \eqref{Eq:ODEy}, we have that
\begin{align*}
	\trace(\mathcal A) = 0 + \trace(\mathbf A).
\end{align*}
Under our assumptions, $\mathcal A$ has only one zero eigenvalue and, hence, the spectrum of $\mathbf A$ is composed by the nonzero eigenvalues of $\mathcal A$. All these eigenvalues have negative real part. 

Although the triplet $(\mathbf A, \mathbf B, \mathbf C)$ in \eqref{Eq:ODEStable} might not be a minimal realisation, it is always possible to obtain a model which is both controllable and observable via its Kalman decomposition \cite{kalman1963mathematical, chen1998linear}. In fact, the command \emph{balreal} of Matlab's Control System Toolbox will obtain the controllable and observable system before obtaining the balanced realisation; hence, is not absolutely necessary to test for these properties separately, when using this software. Thus for stable systems we can perform the model balancing described in Section \ref{Sec:Balanced} in order to obtain a reduced-order model of the form \eqref{Eq:ODEBalRed}. 

Up to now, we had considered that the number of states, $\nums$, of the Markov chain is finite. However, when considering chemical reaction networks, it is possible to use of the Finite State Projection (FSP) method \cite{Munsky2006}, to obtain an ODE set analogous to \eqref{Eq:ODEy} with the most representative, finite number of states. Due to its approximate nature, the set of ODEs obtained via the FSP might not present the properties in \eqref{Eq:Properties}. Hence the change of variables in \eqref{Eq:Trans} would no longer be necessary and balanced model reduction can be applied directly to the FSP-reduced set of ODEs.

Although the lower-dimensional model can be used for obtaining an approximated numerical solution for the probabilities of interest, we would like to remark that one may use the results in \cite{lopez2013exact} to derive closed-form expressions for these probabilities. In the following section, we study some case studies to show the applicability of these methods. 
\section{Case Studies}\label{Sec:Cs}
In this section, we show the derivation and application of reduced order models, through different examples. We will first analyse, in Section \ref{Sec:Monomolec}, one monomolecular reaction and obtain an accurate approximation for the probability of having the conversion of all the molecules from the first species to the second one. Later, in Section \ref{Sec:SMM}, we derive reduced order models capable of approximating a catalytic conversion of a substrate even in cases in which the stochastic Michaelis-Menten cannot yield accurate results \cite{sanft2011legitimacy}. Finally, in Section \ref{Sec:Range} we revisit the catalytic substrate conversion to derive the probability of having ranges of product molecules. In all case studies, we used a 3.2 GHz Quad-Core Intel Xeon computer with 16GB of RAM. Our script was coded in MATLAB\copyright\, R2012b.

\subsection{Monomolecular Reaction Network}\label{Sec:Monomolec}
Let us consider the reversible reaction
\begin{equation}\label{Eq:RRever}
 \ReactionR{\Sb_1}{\Sb_2}{\kf{}}{\kb{}},
\end{equation}
along with the vector composed of species' molecular number $\Sub(\te) := \begin{pmatrix} \sub_1(\te) & \sub_2(\te) \end{pmatrix}^T$. Furthermore, consider an initial number of molecules $\sub(0) = \begin{pmatrix} 300& 0 \end{pmatrix}^T$. Note that in the reaction above the number of molecules remains constant and equal to the initial $300$ molecules. Hence the set of $\SUB $ has $w=301$ elements and may be ordered as follows
\begin{align*}
\SUB = \left \{  \begin{pmatrix}  300 \\ 0   \end{pmatrix},\begin{pmatrix}  299\\1  \end{pmatrix}, \hdots, \begin{pmatrix}  1\\299  \end{pmatrix},\begin{pmatrix}  0 \\ 300  \end{pmatrix} \right \}.
\end{align*}
Now, we are interested in the time-dependent probability of having $300$ molecules of $\Sb_2$, i.e., to be in state $\Sub^{301} = \begin{pmatrix} 0 & 300 \end{pmatrix}^T$. With this formulation, the matrix $\mathcal A$ in \eqref{Eq:ODEy} is shown in Equation \eqref{Eq:DefATrid} on page \pageref{Eq:DefATrid}.
\begin{figure*}[htb!]
\begin{align}
 \mathcal A \iga \begin{pmatrix} -300 \kf{} & \kb{} & 0 &0&\hdots  &0&0 \\
													 \phantom{-}300 \kf{} & -\left (  299 \kf{} + \kb{} \right ) & 2 \kb {} & 0&\hdots & 0&0\\
													 0 & 299 \kf{} &-\left ( 298\kf{} + 2\kb{} \right ) &3\kb{}& \hdots& 0&0\\
													 \vdots&\ddots&\ddots&\ddots&\ddots&\ddots&\vdots\\
													 0&0&0&0&\ddots&-\left ( \kf{} + 299\kb {} \right )&\phantom{-}300\kb{}\\
													 0&0&0&0&\hdots&\kf{}&-300\kb{}
 							\end{pmatrix} \in \mathbb R^{301^2}.\label{Eq:DefATrid}
\end{align}
\end{figure*}
In turn $\Prob(0)$ and $\mathcal C$ are given by
\begin{subequations}\label{Eq:DefPC}
\begin{align}
 \Prob(0) \iga \begin{pmatrix}  1 &0&\hdots&0 \end{pmatrix}^T\in \mathbb R^{301}, \\
 \mathcal C \iga \begin{pmatrix} 0&0&\hdots&1 \end{pmatrix} \in \mathbb{R}^{1\times 301}.
\end{align}
\end{subequations}

With the definitions for $\mathcal A$, $\Prob(0)$, and $\mathcal C$  in \eqref{Eq:DefATrid} and \eqref{Eq:DefPC}, respectively, and by choosing the parameters $\{\kf{}, \kb{} \} = \{150,1\}[\mathrm {s}^{-1}]$, we implemented the model in \eqref{Eq:ODEStable} in Matlab 2012b and obtained its balanced realisation with the command \emph{balreal}. \reffig{Fig:HSV} shows the first $30$ Hankel Singular Values of the balanced realisation's grammian. We observe that the first ten singular values have a large norm in comparison to the rest. By using the command \emph{modred}, we obtained the reduced order model with different number of states; hence, achieving different degrees of approximation. 

\ImagenA{Largest Hankel Singular Values of the balanced realisation of the model of the form in \eqref{Eq:ODEStable}, where $\mathcal A$ and $\mathcal C$ are defined in \eqref{Eq:DefATrid} and (\ref{Eq:DefPC}b), respectively. Additionally $\{\kf{}, \kb{} \} = \{150,1\}[\mathrm {s}^{-1}]$.}{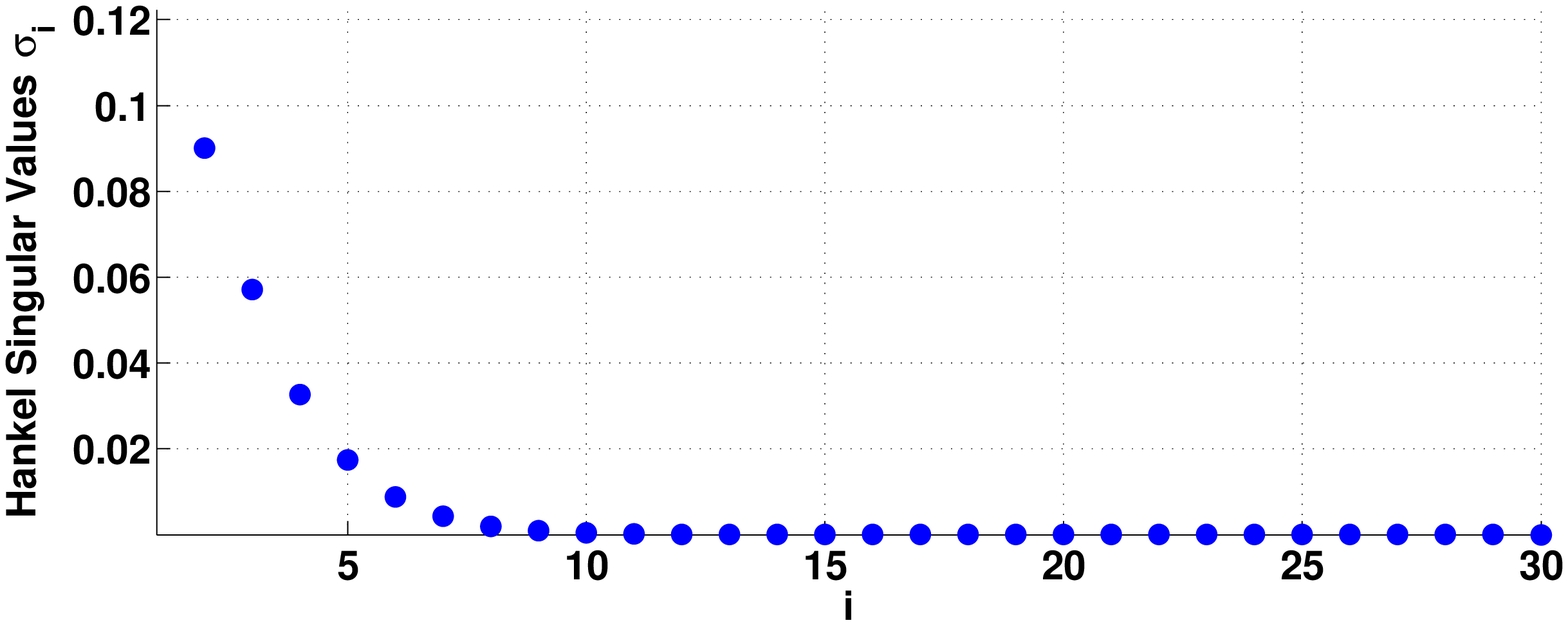}{Fig:HSV}{!p}{0.5}

We depict the impact of the number of states on the error of approximation, in \reffig{Fig:CompReduction}. There, we note that a very coarse approximation is achieved when we try to approximate the full model with $301$ states with a model of only $1$ state (see the lower panel of \reffig{Fig:CompReduction}\textbf{(A)}). In turn, when the reduced order model has $10$ states, the error of approximation is of order $10^{-5}$, as depicted in the lower panel of \reffig{Fig:CompReduction}\textbf{(C)}. Furthermore, if the reduced model has $15$ states, the approximation error might already range in the order of the integration error, as suggested by the irregular fluctuations shown in the lower panel of \reffig{Fig:CompReduction}\textbf{(D)}.

\ImagenA{Output comparison of the full CME and the reduced order model. The upper panels depict of the probability of having all the molecules of $S_1$ converted to $S_2$ by means of the reversible reaction \eqref{Eq:RRever}. The discontinuous line represents this probability as obtained with the full model and the continuous lines with the reduced order model.  In turn, the lower panels show the difference of full model output and that of the reduced order model. The order of the lower-dimensional model for columns \textbf{(A), (B), (C)}, and \textbf{(D)} are $1$, $5$, $10$, and $15$ states, respectively. The parameters used for simulations are as in \reffig{Fig:HSV}.}{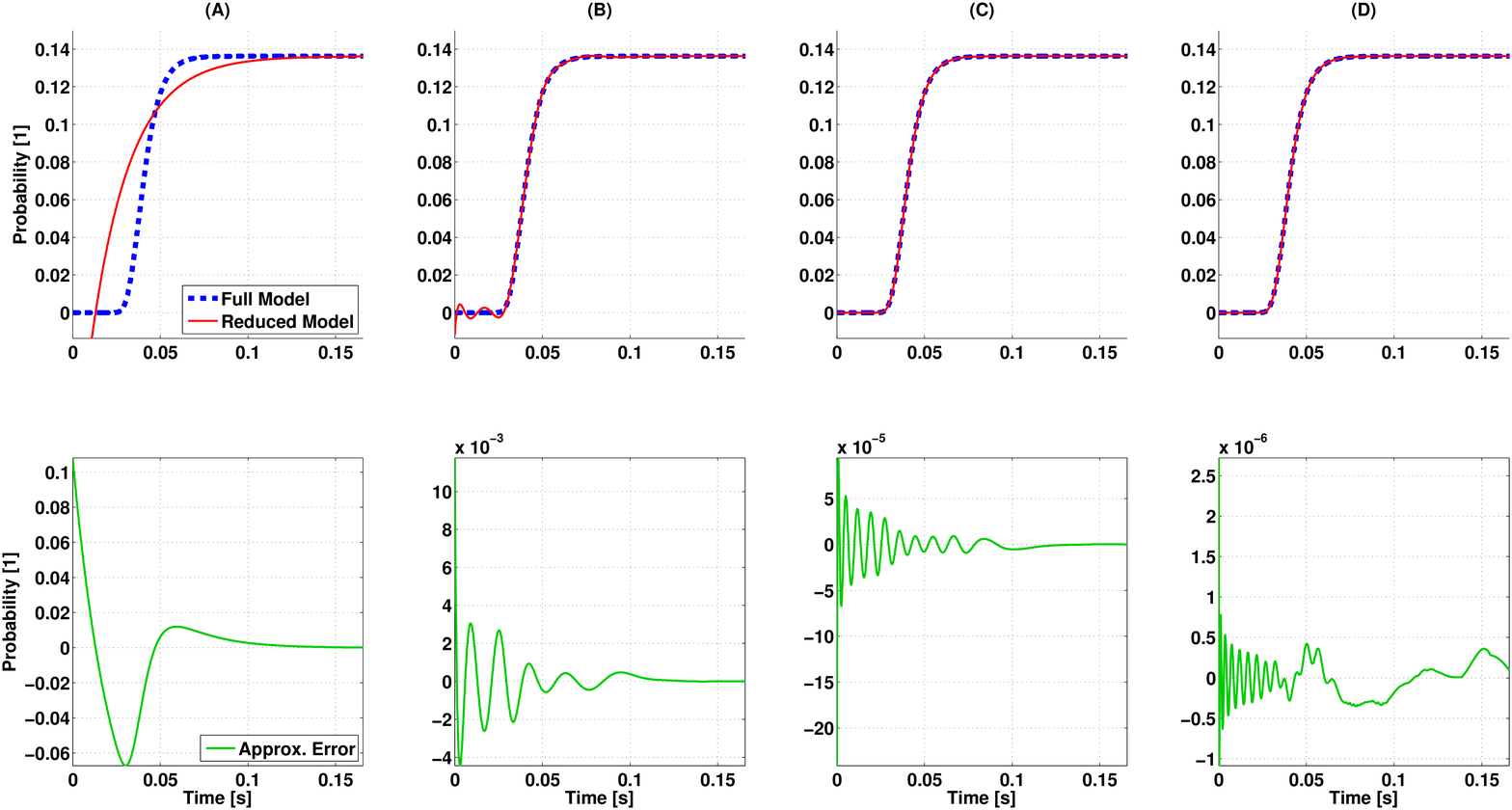}{Fig:CompReduction}{!htb}{1}

To finalise this section, we note that $\mathcal L_2$ gain of the  approximation error is $427.4607 \times 10^{-3},~33.1963\times10^{-3},~587.9172\times10^{-6},$ and $6.0955 \times 10^{-6}$, for the reduced models with $1$, $5$, $10$, and $15$ states, respectively. These bounds were obtained by evaluating Expression \eqref{Eq:ErrorBound}. We note that this is a theoretical bound and does not account for numerical errors during the integration or computation of the Hankel Singular Values.

In the forthcoming section, we obtain reduced order models for a catalytic substrate conversion, and asses the computational burden required to obtain the reduced order model. In addition, we benchmark the time required for simulating the reduced order model against both the computational load required to simulate the full order model and the Stochastic Simulation Algorithm (SSA).

\subsection{Stochastic Michaelis-Menten}\label{Sec:SMM}
In this section, we consider the reaction network 
\begin{align}\label{Eq:CatFull}
		\Reaction{\ReactionR{S+E}{C}{\kf{1}}{\kb{1}}}{P+E}{\kf{2}},
\end{align}
which represents conversion of a substrate, $S$, to a product, $P$, mediated by a catalytic agent, $E$, which binds to the substrate to form the complex $C$. In the deterministic case, it is common practice to approximate the  mass-action-based reaction network in \eqref{Eq:CatFull} via the reaction 
\begin{align}\label{Eq:CatRed}
		\Reaction{S}{P}{v_{MM}(S)},
\end{align}
with nonlinear reaction rate
\begin{subequations}\label{Eq:vMM}
\begin{align}
	v_{MM}([S]) = \frac{v_{max}}{k_m + [S]}[S],
\end{align}
where $[\circ]$ stands for concentration of the argument and
\begin{align}
	v_{max} :\iga \kf 2 [E]_T,\\
	k_m :\iga \frac{\kb 1 + \kf 2}{\kf 1},\\
	[E]_T 	:\iga [E](\te) + [C](\te), 	\quad \textnormal{for any }\te \geq 0. 
\end{align}
\end{subequations}
It has been shown that the dynamics of \eqref{Eq:CatFull} can be reasonably approximated by \eqref{Eq:CatRed} when \cite{segel1989quasi}
\begin{align}\label{Eq:ConMM}
[E]_T << S(0) + k_m.
\end{align}

However, for cases in which the reactions in \eqref{Eq:CatFull} are better described by a stochastic model, it is still possible to represent the dynamics of $S$ and $P$ with a reaction of the form \eqref{Eq:CatRed} by using the propensity
\begin{align}\label{Eq:NlProp}
 a_{MMs}(S) = \frac{v_{max}}{k_m + S}S,
\end{align}
where now $S$ represents the number of molecules of the substrate and $v_{max}$, $k_{m}$, and $E_T$ are as in \eqref{Eq:vMM}, while the kinetic constants are those of the stochastic model. This representation is valid under the condition \eqref{Eq:ConMM}, as is also true for the deterministic case\cite{sanft2011legitimacy}. Additionally, the case when 
\begin{align}\label{Eq:ConMM2}
\kb 1 >> \kf 2, 
\end{align}
in which a nonlinear propensity function can be used to represent \eqref{Eq:CatFull} via \eqref{Eq:CatRed} was considered in \cite{sanft2011legitimacy}. There, the condition $\kb 1 >> \kf 2$ induces a time-scale separation which is further used to abridge the reaction network \eqref{Eq:CatFull}. In the following, we will refer to the representation of the reaction network \eqref{Eq:CatFull} via \eqref{Eq:CatRed} with the propensity in \eqref{Eq:NlProp} or that in \cite{sanft2011legitimacy} as \textit{stochastic Michaelis-Menten} representation.

We now obtain a reduced order model that approximates the probability of being in selected states of the underlying Markov chain. We derive this reduced model by means of the procedure described in Section \ref{Sec:OrderRed}. In contrast to the approaches in \cite{sanft2011legitimacy} and references therein, we do not assume any particular relation among the parameters and initial conditions, so our methodology is more widely applicable.

Another difference from the approaches in \cite{sanft2011legitimacy} is that they prove the applicability of SSA algorithms with the stochastic Michaelis-Menten propensity. In contrast, we derive a dynamical system that approximates the solution of the CME with an \emph{a priori} error bound given by \eqref{Eq:ErrorBound}. We recall that in the limit, the probability distribution obtained from the SSA trajectories will converge to the solution of the CME. However, depending on the kinetic parameters and network analysed, the SSA might require large computational times to provide results with the desired accuracy.

It is noteworthy that even when \eqref{Eq:CatFull} can be represented by \eqref{Eq:CatRed}, one can still obtain a reduced model via the balanced model reduction described in the Section \ref{Sec:OrderRed}, as we do not assume any relationship among the parameters and initial conditions. 

To exemplify the concepts above, we depict in \reffig{Fig:ReductionMM} a comparison of 
\begin{inparaenum}[\bfseries (A)]
\item the solution of the CME of \eqref{Eq:CatFull} with propensities shown in Table \ref{Tb:Prope}; 
\item the solution of the CME of the stochastic Michaelis-Menten in \eqref{Eq:CatRed} with the nonlinear propensity in \eqref{Eq:NlProp}; and
\item the solution of the reduced model described in Section \ref{Sec:OrderRed} for the last state of the Markov chain, which represents total conversion of the substrate to product.
\end{inparaenum}
The parameters used are  $\{\kf 1, \kf 2, \kb 1 \} = \{1~[\left ( \mathrm{molecules ~ s}  \right )^{-1}],1~[\mathrm s ^{-1}],1~[\mathrm s ^{-1}]\}$ and 10 initial molecules of substrate. With these parameters, the method proposed in \cite{sanft2011legitimacy} would not yield accurate results, as the condition \eqref{Eq:ConMM2} is not fulfilled ($\kb 1 = \kf 2$).  

The only difference between the upper and lower panels in \reffig{Fig:ReductionMM} is the number of initial molecules considered for the enzyme. In the upper panel we considered 1 molecule of the enzyme, hence condition \eqref{Eq:ConMM} is fulfilled, and the stochastic Michaelis-Menten representation may be used to approximate the full model. Moreover, one can use the stochastic Michaelis-Menten to derive a reduced model via balanced realisation, as compared in the upper panel of the upper panel of \reffig{Fig:ReductionMM} \textbf{(C)}. There, we approximated the stochastic Michaelis-Menten model of $11$ states with a reduced order model with $6$ states; the $\mathcal L_2$ gain of the approximation error is less than $0.2807 \times 10^{-3}$ as given by \eqref{Eq:ErrorBound}. In contrast, when we consider $10$ molecules of enzyme initially, condition \eqref{Eq:ConMM} is violated (as $E(0) = S(0)$) and the stochastic Michaelis-Menten model does not reproduce the dynamics of the full reaction network in \eqref{Eq:CatFull}, as depicted in the lower panels of \reffig{Fig:ReductionMM}. We note, however, that we can still obtain a reduced model via balanced realisation that accurately approximates the dynamics of the full model (cfr. \reffig{Fig:ReductionMM}\textbf{(C)} lower). There we approximated the full model with $66$ states by a reduced model of $6$ states, whose approximation error $\mathcal L_2$ gain is less that $0.21547\times10^{-3}$.

\ImagenA{Validity of the Michaelis-Menten propensity as an approximation of a catalytic substrate conversion. Column \textbf{(A)} shows the simulation of the CME associated to \eqref{Eq:CatFull}, where each tread represents the probability of being in every state of the Markov chain; in turn, column \textbf{(B)} shows the solution of the CME of \eqref{Eq:CatRed} by using the nonlinear propensity function \eqref{Eq:NlProp}; whereas, column \textbf{(C)} shows the probability of being in the last state of the Markov chain, which represents total conversion of the substrate to the product. This probability is obtained via the CME of the full stochastic model, by the CME of reduced Michaelis-Menten stochastic model, and by the approximated model to the CME via balanced realisation. The parameters used for obtaining the numerical solution are $\{\kf 1, \kf 2, \kb 1 \} = \{1,1,1\}$ and $10$ initial molecules of substrate. The only difference between the upper and lower panels is the number of enzymes considered: upper panels $1$ molecule, whereas the lower panels, $10$ molecules. Note that in the lower panel the stochastic Michaelis-Menten approximation is not valid, but the approximation via the balanced model truncation is close to the full model.}{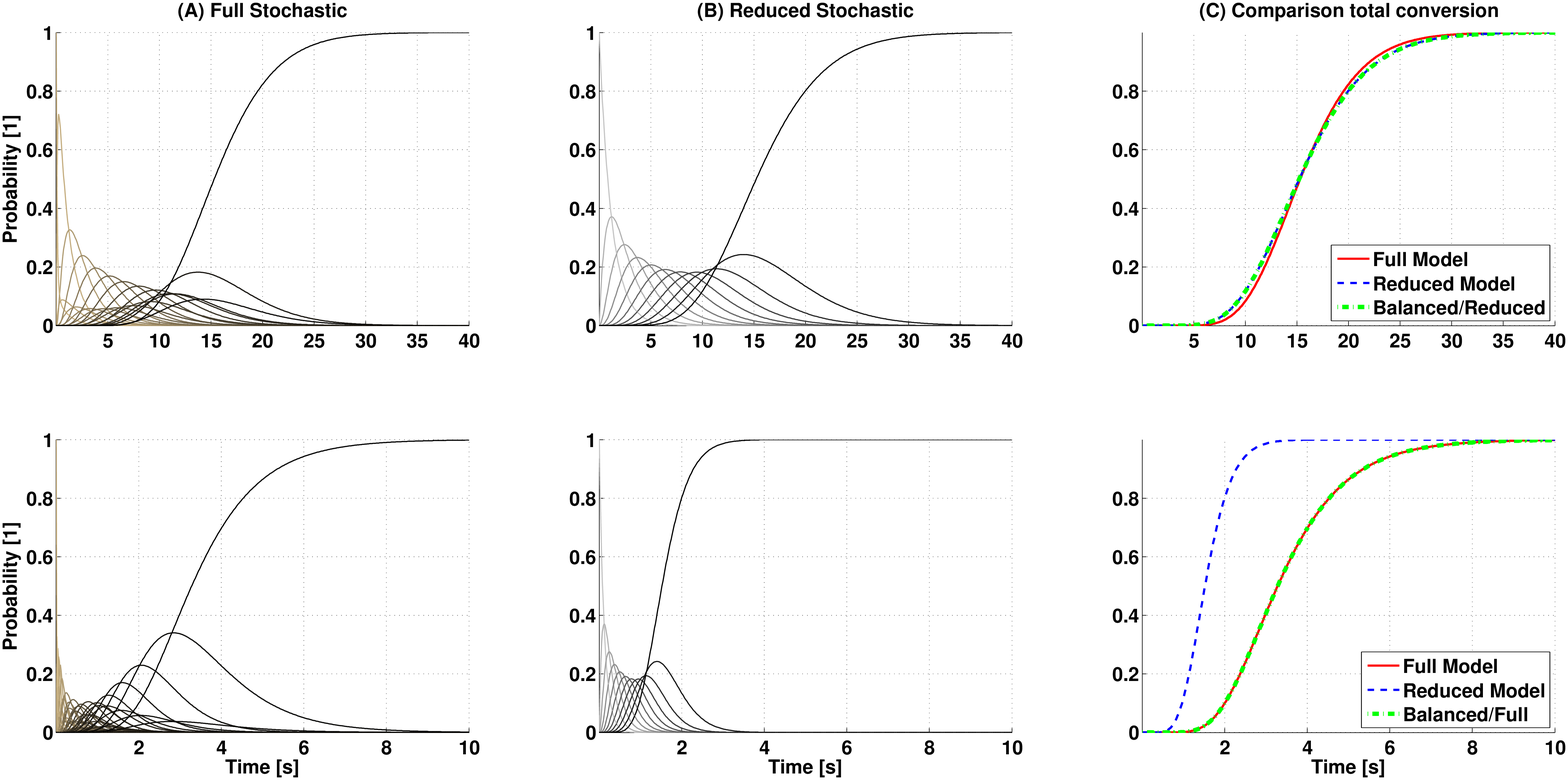}{Fig:ReductionMM}{t!}{1}

Now we focus on the time required to simulate the CME and the time required to simulate the reduced order model. To compute the latter, we need to apply some state transformations to the CME \eqref{Eq:ODEy} to derive a balanced realisation that can be further truncated. 
%
%
Once the reduced model is obtained, the time required for its numerical solution is significantly smaller compared to the time required for the numerical solution of the full CME. To illustrate this reduction on the computational time, we obtained the CME of the reaction network \eqref{Eq:CatFull} with an equal number of molecules for the substrate and enzyme and zero molecules for the rest of the species, in the initial state; later, we obtained the reduced order model via balanced realisation, which represents the state of total conversion of the substrate to the product; and compared the time required for obtaining the numerical solution of the full CME ($\te_{\textrm{CME}}$) and the reduced model ($\te_{\textrm{red}}$) by the expression 
\begin{align}\label{Eq:Per}
 	\eta = \log_{10} \left ( \frac{\te_{\mathrm{CME}} - \te_{\mathrm{red}}}{\te_{\mathrm{red}}}\right ).
\end{align}
We depict the results of this assessment in \reffig{Fig:Performance}. There we observe that as the number of molecules for $E$ and $S$ in the initial state increase, the savings on the computational time required to obtain the numerical solution of the lower-order model also increases. We note that for the comparison in \eqref{Eq:Per} we did not account for the time required to obtain the reduced order model. 

\ImagenA{Computational time overhead, as given by \eqref{Eq:Per}, required to solve the full CME (diamonds) and to perform $10^3$ SSA runs (squares) as compared to the computational time required in seconds to simulate the reduced order model, as the initial number of molecules for $E$ and $S$ vary from $5$ to $100$. The parameters values used for simulation are identical to those of \reffig{Fig:ReductionMM}.}{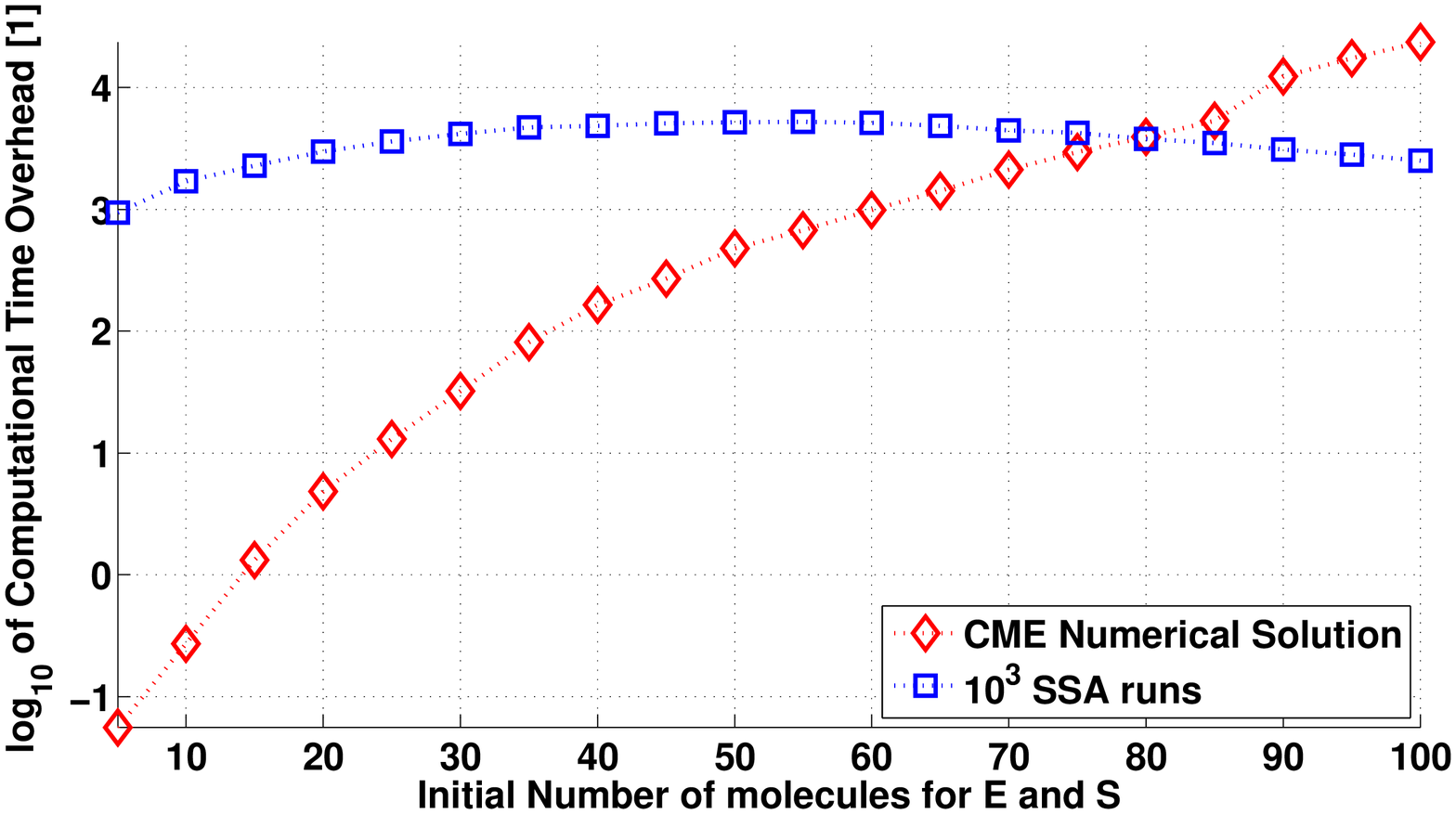}{Fig:Performance}{p!}{0.5}

To finalise this section, we compare the computational time required by 
\begin{inparaenum}[\itshape i\upshape)]
\item the derivation of the reduced model via balanced realisation plus the simulation of the reduced model; and
\item the time required by the FSP \cite{Munsky2006} for each time point. 
\end{inparaenum}  %
We note that the FSP obtains an approximated probability vector with a desired error bound ($\varepsilon$) for \textit{one specific time point}; hence, if one is interested in the transient response of the probability distribution, one has to implement such an algorithm for every time step of interest. In contrast, once one obtains the reduced model via balanced realisation it is possible to use the lower-dimensional system for any number of time points. These results are summarised in \reffig{Fig:FSPvsBalReal}, where the panels $\mathbf{(A)}$, $\mathbf{(B)}$, and $\mathbf{(C)}$ consider $10$, $30$, and $50$ initial molecules for $E$ and $S$, respectively, and zero molecules for the rest of the species. The remaining parameter values are identical to those of \reffig{Fig:ReductionMM}. We note that for the FSP the 1-norm of the error bound is less than a predefined $\varepsilon$ for the specific time points of interest (discrete signal), whereas the $\mathcal L_2$ gain of the approximation error (continuous signal), obtained with the reduced model via balanced realisation, satisfies the bound given by \eqref{Eq:ErrorBound}. As the nature of both error signals is different, is difficult to perform a fair comparison of the methods' accuracy. In the forthcoming section, we obtain a reduced order model that approximates the probability of having a certain range of $P$ molecules.

\ImagenA{Comparison of the computational time required to obtain the reduced order model via balanced realisation (filled circle) and to obtain the approximative model via the FSP method (empty markers), with different, predefined error bounds ($\varepsilon$). The reaction network analysed is \eqref{Eq:CatFull}. The parameters used for simulation are those of \reffig{Fig:ReductionMM}. Panels $\mathbf{(A)}$, $\mathbf{(B)}$, and $\mathbf{(C)}$ consider $10$, $30$, and $50$ initial molecules for $E$ and $S$ and zero molecules for the rest of the species, respectively.}{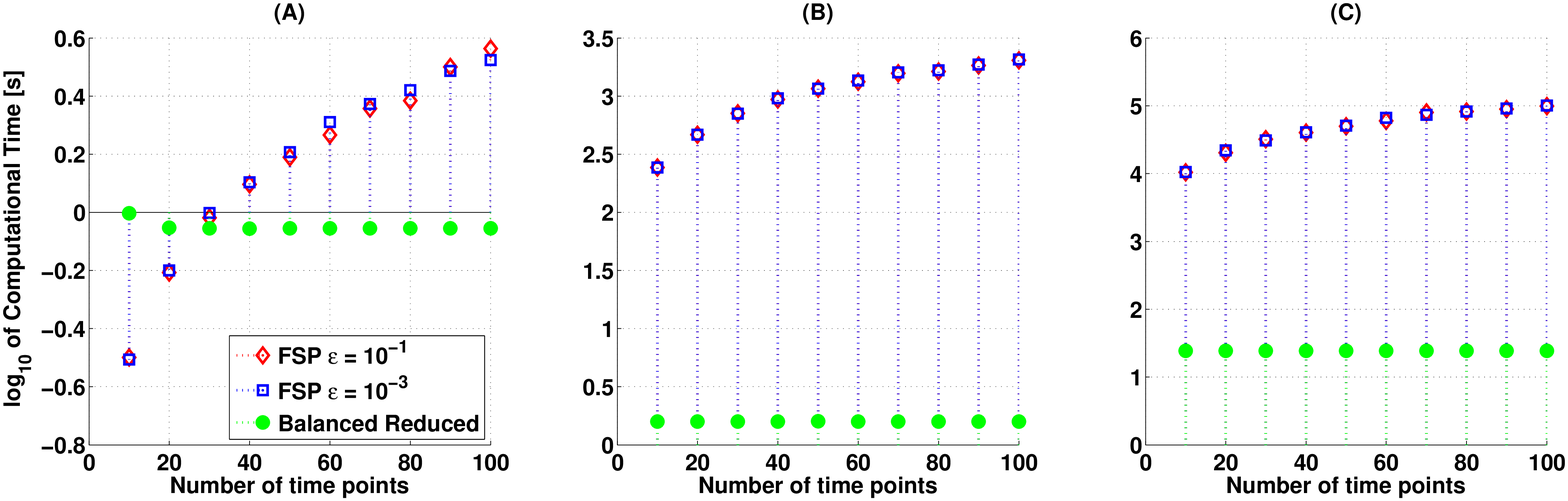}{Fig:FSPvsBalReal}{p!}{1}

\subsection{Probability for Ranges of Molecules Counts}\label{Sec:Range}
Up to now, we have obtained reduced models that 	approximate the probability of being in one state of the Markov chain. In this section, we revisit the reaction network in \eqref{Eq:CatFull} by obtaining the probability of having a certain number of molecules within predefined ranges. Here we consider the following parameter definitions: $\{\kf 1, \kf 2, \kb 1 \} = \{1,1,1\}$, $100$ initial molecules of substrate, $100$ initial molecules of enzyme, and zero initial molecules for the rest of the species. By denoting the number of $P$ molecules with $\con_P$, we can formulate our problem as approximating the following probabilities
\begin{align}\label{Eq:ypart}
 \mathbf{y}(\te) = \begin{pmatrix}  \Pr(0 \leq \con_P(\te) \leq 30) \\ \Pr(31 < \con_P(\te) \leq 70)  \\  \Pr(71 < \con_P(\te) \leq 100)
 \end{pmatrix}.
\end{align}
To derive the CME, one needs to obtain and label all the possible combinations of species molecular counts $\Sub^i$ and organise them in the set $\SUB$ in \eqref{Eq:SetSpace}. Then we have to evaluate the infinitesimal generator $\mathcal A$ as in \eqref{Eq:PropMat} with the corresponding reaction propensities of \eqref{Eq:CatFull} (see Table \ref{Tb:Prope}). To obtain an expression for $\mathbf y(\te)$, we need to define the matrix $\mathcal C$ in \eqref{Eq:Defye} so that the product of the first row of $\mathcal C$ by the vector $\Prob(\te)$ yield the sum of the probability of all the states $\Sub^i$ such that $\con_P$ is within the range $[0,30]$. The next two rows of $\mathcal C$ are defined likewise, but accounting for the ranges $\con_P$ described in the second and third entries of \eqref{Eq:ypart}. The CME for this system, parameters, and initial number of molecules has $5151$ states. By applying the model reduction technique in Section \ref{Sec:OrderRed}, we can approximate the probabilities in \eqref{Eq:ypart} by a dynamical system with $16$ states, whose output is depicted in \reffig{Fig:Intervals}. The $\mathcal L_2$ gain of the approximation error is less than $6.384\times 10^{-3}$, as estimated by \eqref{Eq:ErrorBound}.

\ImagenA{Marginal probability distributions for the reaction network \eqref{Eq:CatFull}. Probability of having a molecular count of $P$ within a certain range, as obtained with the reduced order model. The parameters used for simulation are $\{\kf 1, \kf 2, \kb 1 \} = \{1,1,1\}$, $100$ initial molecules of substrate, $100$ initial molecules of enzyme, and zero initial molecules for the remaining species.}{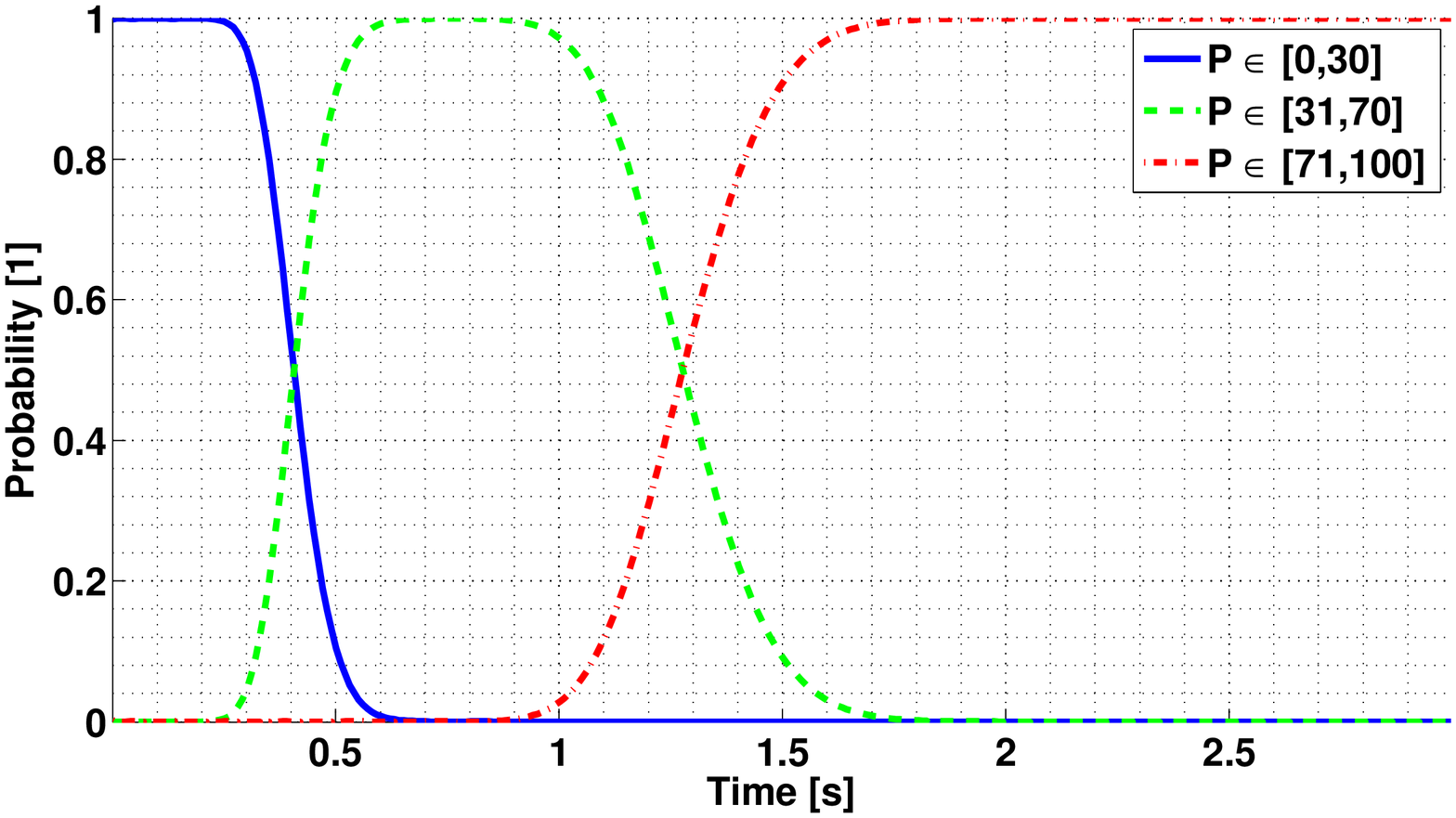}{Fig:Intervals}{p!}{0.5}
\section{Concluding Remarks}\label{Sec:Con}
In this paper we addressed the order reduction of the infinitesimal generator of a homogeneous, continuous-time, finite and discrete state-space Markov Chain via the reduction of its balanced realisation. Although the application range of these dynamical systems is broad, here we focus on its use on stochastic chemical reaction networks, without loss of generality. In this context, the infinitesimal generator of the Markov chain that describes the probability of having a particular species molecular count is a large set of ODEs. 

To reduce the order of the infinitesimal generator of a Markov chain, we used an alternative coordinate system to represent the Chemical Master Equation (CME). This representation, denoted as Lyapunov balanced realisation, has interesting property that the states are organized in decreasing order according to the probabilities of interest. 
Hence, an accurate approximation can be obtained, for example, by neglecting the last states of the Lyapunov balanced model, as discussed in Section \ref{Sec:OrderRed}. Although one may focus on particular states of the Markov chain, it is also possible to account for marginal probability distributions such as in the case study in Section \ref{Sec:Range}, or even mean values, by properly defining the matrix $\mathcal C$ in \eqref{Eq:Defye}.

In many cases, only selected states of the Markov chain might be of practical relevance. For instance, this is the case when facing limited or inexact measurement data, or when only a few states are relevant for downstream signalling in biochemical reactions. Also, in imaging analysis of chemical reaction networks, obtaining the exact count of intracellular protein reporters might be challenging due to limited resolution. Hence, the validation of the mathematical model that describes the process under observation should yield the probability of having an specific range of molecules count of the observed species. We presented this procedure in Section \ref{Sec:Range}, for a very simple reaction network. Even in such a simple case, the associated Markov chain presented approximately 5000 distinct states of the system. This highlights how simulation of a system, even in the simplest cases, might imply a computationally intensive task. To alleviate such a burden, the model reduction via balanced realisation used in this paper yields lower-dimensional ODE sets, 
whose numerical solution might be several orders of magnitude quicker than the numerical solution of the original CME. Moreover, the method used to derive the lower dimensional model provides an upper bound on the approximation error, depending on the number of states neglected to derive the approximation.

Of note, the processes required for deriving the reduced order model itself might take longer computations times compared to the mere simulation of the CME. 
Nevertheless, the numerical solution of the reduced model might be obtained orders of magnitude faster, depending on the number of molecules of the system, as shown in \reffig{Fig:Performance}. Hence, there will be real savings on the computational time when the reduced model is repeatedly utilised, for instance when adopting different initial probability distributions. We would like to stress that to obtain a reduced order model, we have to fix kinetic parameters and to define which are the states of interest. Should we require to modify either of them, a new reduced model has to be derived. Likewise, all methods that require computational calculations, such as the FSP, SSA, and numerical solution of the CME will require numeric values for the parameters and, moreover, specific numerical values for the initial probability distribution. When either of them are modified, a new numerical solution has to be obtained. Additionally, the reduction and simulation of the CME might be orders of magnitude faster than the application of the FSP method, as suggested by the example analysed in Section \ref{Sec:SMM}.

Another possible use for the reduced model is to derive closed-form expressions of its solution (see \cite{lopez2013exact}, for instance), thereby avoiding the need for numerical solution of the reduced ODE set. When the number of states of the Markov chain to reduce is so large that using only one computer is unfeasible, we suggest the use of parallel algorithms to obtain the model reduction by truncation (see e.g. \cite{badia2006parallel,benner2005model})

It is important to note that the reduced order model might lack some properties of the full model. For instance, the infinitesimal generator of the Markov chains studied here describes a positive system: the value of the probabilities will be always positive. However, the reduced order model obtained by truncation used in this paper will not, in general, preserve such a property. This implies that if most of the states of the balanced realisation are neglected to obtain the reduced model, there is a risk of having small, negative values for the approximated probabilities. An example of such phenomenon can be observed on the upper panels of \reffig{Fig:CompReduction}\textbf{(A,B)}. This suggests the existence of a trade-off on the order and the accuracy of the reduced-order model. As a rule of thumb, a good approximation can be obtained by neglecting those states associated to Hankel Singular Values which are three orders of magnitudes smaller than the largest one. If the possibility of small, negative values for the probability cannot be afforded for the application of the reduced order model, there are other model order reduction methods that preserve the positivity of the original model, such as the recent works \cite{feng2010internal,li2011positivity,grussler2012symmetry}. However, it is equally important to note that these approaches are not generally applicable; are more time consuming; and have larger error bounds.


Along this paper, we have considered that the set $\SUB$ in \eqref{Eq:SetSpace} has all the possible states of the Markov chain under consideration. However, when the number of states is prohibitively large, it is possible to consider a truncation of the set $\SUB$; thereby, obtaining smaller Master Equations. This truncation has two implications: 
\begin{inparaenum}[\itshape i\upshape)]
\item the Master Equation derived from the truncated $\SUB$ will not capture the full probability density function of the Markov chain, but will only focus on the probability of being in those states of interest as characterised in \cite{Munsky2006}; and 
\item the set of ODEs arising from the truncated $\SUB$ will not have the properties in \eqref{Eq:Properties}.
\end{inparaenum} 
Hence the change of coordinates in \eqref{Eq:Trans} is not necessary, and balanced model reduction can be applied directly to the set of ODEs obtained from the truncated $\SUB$. This, in turn, implies that those methods that depend on the truncation of the set $\SUB$ to derive approximated probability distributions, such as \cite{Munsky2006}, do not antagonise with the model reduction via balanced realisation used in this paper, as both approaches can be complementary.
\appendix
\section{Bound on the Approximation Error of the Model Reduction via Balanced Realisation}\label{Sec:App}
Here, we provide some definitions and the derivation of the approximation error bound \eqref{Eq:ErrorBound} that arises from the model reduction via balanced realisation described in Section \ref{Sec:OrderRed}. The material of this section is based on the Refs. \cite{Zhou1996robust,Skogestad2007}.

First, to asses the size of the error of approximation, let us define the $\mathcal L_2$ norm of a real, time-dependent vector $\mathbf u(\te)$ as
\begin{align*}
	||\mathbf u||_{\mathcal L_2} := \int_0^\tau \mathbf u^T(\te)\mathbf u(\te) \mathrm d t.
\end{align*}
When $\tau < \infty$, one obtains the norm of the truncated signal $\mathbf u (\te)$. To increase readability, we will not explicitly show the upper limit of integration in the norm's subscript. 

Now, in the frequency domain, the linear ODE \eqref{Eq:LinODE} becomes the following algebraic equation
\begin{equation}\label{Eq:YLap}
 \mathbf{Y}(\xi) = \mathbf G(\xi) \mathbf{U}(\xi),
\end{equation}
where $\xi$ is the complex frequency variable that arises from the Laplace transform of \eqref{Eq:LinODE}, and 
\begin{equation*}
 \mathbf{G}(\xi) : = \mathbf{D} + \mathbf{C} \left ( \xi \mathbf {I - A} \right )^{-1} \mathbf{B}.
\end{equation*}
The complex matrix $\mathbf {G}(\xi)$ is denoted as the \emph{transfer function} of the system \eqref{Eq:LinODE} and characterises its input-output behaviour. The $\mathcal H_{\infty}$ norm of the complex matrix $\mathbf {G}(\xi)$ is defined as
\begin{equation*}
 ||\mathbf {G}||_{\mathcal H_{\infty}} := \mathrm{sup}_{\mathrm{Re}(\xi)>0} \sqrt{\bar{\lambda} \left (\mathbf {G}^*(\xi) \mathbf {G}(\xi)\right )}. 
\end{equation*}
Here $\bar{\lambda} \left ( \circ \right )$ denotes the largest eigenvalue of the argument. In turn, the $\mathcal H_{2}$ norm of $\mathbf {G}(\xi)$, for analytic matrices on the open right half-plane, is
\begin{equation*}
||\mathbf U||_{\mathcal H_{2}} = \frac{1}{2\pi}\int_{-\infty}^{\infty}  \mathrm{Trace} \left ( \mathbf{Y(\mathnormal j\omega)^* Y(\mathnormal j \omega)}   \right )  \mathrm d \omega.
\end{equation*}
When $\mathbf{U}(\xi)$ in \eqref{Eq:YLap} belongs to the Banach space endowed of the norm $\mathcal H_2$, Theorem 4.4 in \cite{Zhou1996robust} states that 
\begin{align}\label{Eq:GainH}
 	||\mathbf Y||_{\mathcal H_{2}}^2 = ||\mathbf {G}||_{\mathcal H_{\infty}}^2 ||\mathbf U||_{\mathcal H_{2}}^2.
\end{align}

In order to relate the frequency-domain norms with the time-domain norms, we note that the Laplace transform used to obtain the transfer function of \eqref{Eq:LinODE} is an isomeric isomorphism between the $\mathcal H_2$ space in the frequency-domain and the $\mathcal L_2$ space in the time-domain. Thus, from \eqref{Eq:GainH}, we can infer that 
\begin{align}\label{Eq:GainL}
 	||\mathbf y||_{\mathcal L_{2}}^2 = ||\mathbf {G}||_{\mathcal H_{\infty}}^2 ||\mathbf u||_{\mathcal L_{2}}^2.
\end{align}

Now, we are ready to present the error bound due to the model-order reduction as presented in \cite[Th. 11.1]{Skogestad2007}
\begin{Theorem}
Let $\mathbf{G}(\xi)$ be a stable rational transfer function with Hankel singular values $\sigma_1 \geq \sigma_2 \geq ...\geq \sigma_\nums$ and let $\mathbf{G}_{\mathrm{red}}(\xi)$ be obtained by truncating or residualising the balanced realization of $\mathbf{G}(\xi)$ to the first $k$ states. Then 
\begin{equation}
	||\mathbf {G} - \mathbf{G}_{\mathrm{red}}||_{\mathcal H_{\infty}} \leq 2  \sum_{i=k+1}^{\nums} \sigma_{i}.
\end{equation}
\end{Theorem}
Hence, the relationship in \eqref{Eq:GainL} implies
\begin{align}\tag{\ref{Eq:ErrorBound}}
 	\frac{\abs{\abs{ \mathbf y  -  \mathbf y_{\textrm{red}}}}_{\mathcal L_2}}{\abs{\abs{\mathbf u}}_{\mathcal L_2}} \leq 2  \sum_{i=k+1}^{\nums} \sigma_{i}.
\end{align}

\end{document}